\documentclass[12pt]{amsart}

\usepackage{latexsym}
\usepackage{amsmath}

\newtheorem{theorem}{Theorem}
\newtheorem{conjecture}[theorem]{Conjecture}
\newtheorem{corollary}[theorem]{Corollary}
\newtheorem{lemma}[theorem]{Lemma}

\newtheorem{proposition}[theorem]{Proposition}
\newtheorem{remark}[theorem]{Remark}
\newtheorem{question}[theorem]{Question}

\DeclareMathOperator{\Gdim}{Gdim}
\DeclareMathOperator{\cof}{cof}

\begin{document}
\title{On the Goldie Dimension of Hereditary Rings and Modules}

\author{H.Q. Dinh}
\address{Department of Mathematical Sciences, Kent State
University - Trumbull, Warren, OH 44483, USA}
\email{hdinh@kent.edu}

\author{P.A. Guil Asensio}
\address{Departamento de Matem\'{a}ticas Universidad de
Murcia, 30100
Espinardo, Murcia, Spain}
\email{paguil@um.es}

\author{S.R. L\'opez-Permouth}
\address{Department of Mathematics, Ohio University,
Athens, OH 45701, USA}
\email{slopez@math.ohiou.edu}

\thanks{The second author  has been partially supported
by the DGI (BFM2003-07569-C02-01, Spain) and by the Fundaci\'on
S\'eneca (PI-76/00515/FS/01)}



\begin{abstract}
We find a bound for the Goldie dimension of hereditary modules in
terms of the cardinality of the generator sets of its
quasi-injective hull. Several consequences are deduced. In
particular, it is shown that every right hereditary module with
countably generated quasi-injective hull is noetherian. Or that
every right hereditary ring with finitely generated injective hull
is artinian, thus answering a long standing open question posed by
Dung, G\'omez Pardo and Wisbauer.
\end{abstract}

\maketitle

\section{Introduction.}

Let $R$ be an associative ring with identity. The (infinite)
Goldie (or uniform) dimension of a right $R$-module $M$, $\Gdim(M)$,
is defined to be the supremum of all cardinal numbers $\aleph $
such that there exists a direct sum $\oplus _{I}M_{i}\subseteq M$
of non-zero submodules
of $M$ with $%
\left\vert I\right\vert =\aleph \allowbreak $ (see e.g.
\cite{DF}). This definition is based on lattice-theoretical
properties of the set of submodules of $M$. Therefore, it can
also be stated in terms of Lattice Theory, specifically in
modular lattices (see \cite{SS}). This gives the context of an
interesting common framework where it is possible to include, for
example, the study of eigenvalues of Hermitian compact operators
or the singular values of compact operators (see \cite{ISS}). It
is worth mentioning that these theories give examples of lattices
of infinite Goldie dimension (concretely, countable) that are not
covered by the classical notion of finite Goldie dimension.

A main question when dealing with infinite Goldie dimensions of
modules (or modular lattices) is whether a module of Goldie
dimension $\aleph$ must necessarily contain a direct sum of
$\aleph $ non-zero submodules. In this case it is said that
$\Gdim(M)$ is \emph{attained in} $M$. The main result of \cite{DF}
(see also \cite{SS}) asserts that $\Gdim(M)$ is attained whenever
it is not a weakly inaccessible cardinal. Moreover, some examples
are given showing that for weakly inaccessible cardinals, Goldie
dimensions may be not attained. These examples are based on
constructions by \"{E}rdos and Tarski for Boolean Algebras
\cite{ET}. However, let us point out that the existence of a
weakly inaccessible cardinal would imply the consistency of ZFC
(Zermelo-Fraenkel Set Theory with the Axiom of Choice), even
assuming the Generalized Continuum Hypothesis (see, e.g.
\cite{HJ},\cite{R}). Therefore, due to G\"{o}del Incompleteness
Theorems, the existence of a weakly inaccessible cardinal cannot
be proved in ZFC. Actually, the existence of this type of
cardinals is an independent statement in ZFC.

Using certain combinatorial methods in Set Theory, commonly
refereed to by many Ring theorists as ''Tarski's Lemma''
\cite{T}, it
was proved in \cite{GG} that if $M$ is a non-singular injective (or, more generally,
a quasi-continuous) module which contains an essential finitely
generated submodule and if $\Gdim(M)=\aleph$ is attained and
infinite, then there exists a quotient of $M$ that contains an
infinite direct sum of $\aleph ^{+}$ non-zero submodules, where
$\aleph ^{+}$ stands for the successor cardinal of $\aleph $.
This technical result was then applied to assure the
finiteness of the Goldie dimension of $M$ under various assumptions. Consider, for
example, when every quotient of $M$ has countable Goldie
dimension
\cite[Corollary 2.4]{GG}.

In this paper, we continue this line of reseach by studying the
Goldie dimension of non-singular modules over hereditary rings. We
first observe that the Goldie dimension of a projective module $P$
is closely related to the minimal cardinality of the
generator sets of its submodules. This relation is essential
for proving our main results. First we show that if $P$ is a finitely
generated hereditary module whose quasi-injective hull is
$\aleph$-generated for some infinite cardinal number $\aleph$,
then $M$ cannot contain an independent family of $\aleph$ non-zero
submodules. When $\aleph$ is not an inaccessible cardinal this means that $\Gdim(P)<\aleph$.
Several consequences are derived from
them. Among them, we show that every right hereditary module such
that the quasi-injective hulls of its cyclic submodules are
countably generated is a direct sum of noetherian modules. In
particular, we get that every right hereditary ring with countably
generated injective hull is right noetherian, thus extending
results in \cite{GDW} and \cite{GG1}. Possibly the most
interesting consequence is Corollary \ref{GDW}, where it is shown
that every right hereditary ring with finitely generated injective
envelope is right artinian. This corollary answers an open
question posed in \cite{Du} and \cite[Remark (b), p. 1033]{GDW}
that has remained open for fifteen years and whose motivation
sterns back to Osofsky's work on hypercyclic rings as well as to
an old characterization of (two-sided) hereditary Artinian QF-3
rings given in \cite[Theorem 3.2]{CR}. See \cite[Corollary 6]{GDW}
and \cite[Corollary 1.10]{GG1} for some partial answers to this
problem.

In Section 3 we extend our results to modules $M$ over right
hereditary rings such that the (quasi-)injective hull of $M$ is
$\aleph$-presented for some cardinal number $\aleph$ of cofinality
$\omega$ (see Theorem \ref{singular}). In particular we obtain
that any countably presented (quasi-)injective module over such
rings is a direct sum of uniform modules . Let us point out that,
in order to get this theorem, we need to drop the nonsingularity
hypothesis in $M$, a critical condition in this kind of results
(see \cite{O1,O2} for a deep discussion about this question). On
the other hand any cardinal number of cofinality $\omega $ is
(both weakly and strongly) accessible and, therefore, Goldie
dimension is always attained for these cardinals (see \cite{DF}).

We close the paper by discussing in Section 4 possible ways to
extend our techniques to a more general question. Namely, whether
any finitely generated module satisfying that every quotient is
injective is a direct sum of uniforms (Conjecture \ref{conj}). We
also show that our results give partial positive answers to this
other conjecture.

Throughout this paper all rings will be associative and with
identity. Mod-$R$ will stand for the category of right
$R$-modules. Given a cardinal number $\aleph$, we will denote by
$\aleph^+$, the successor cardinal of $\aleph$, i.e., the smallest
cardinal number that is strictly greater than $\aleph$.

We refer to \cite{AF,D,HJ,W} for any notion used in the text but
not defined herein.

\section{Main Results.}

Let $\aleph $ be cardinal number. A right $R$-module $M$
is called $\aleph $%
-generated if there exists an epimorphism $p:R^{\left( \aleph
\right) }\rightarrow M$. And $M$ is said to be $\aleph $-presented
if the kernel of this epimorphism is also $\aleph $-generated.
Given a module $M$, we are going to denote by $Add[M]$ the full
subcategory of Mod-$R$ consisting on direct summands of direct
sums of copies of $M$. The following easy lemma points out the
close relation that there exists between the cardinality of the
generator sets of projective modules and their decompositions as
direct sums of non-zero submodules. Since any projective module is
an element of $Add[R_{R}]$, we state the lemma in this more
general form.

\begin{lemma}
\label{generators} \label{L1}Let $\aleph$ be an infinite
cardinal number, $M$, an $\aleph$-generated right
$R$-module and $\aleph^{\prime}$, a
cardinal number with $\aleph^{\prime}>\aleph$. Given an $N\in
Add[M]$, the following conditions are equivalent:

\begin{enumerate}
\item $N$ can be decomposed as a direct sum of $\aleph^{\prime}$
non-zero direct summands

\item Every generator set of $N$ has cardinality at least
$\aleph^{\prime}$.
\end{enumerate}

In particular, given an uncountable cardinal number $c$, a
projective module $P$ is the direct sum of $c$ non-zero
submodules if and only if every generator set of $P$ has
cardinality at least $c$.
\end{lemma}

\begin{proof}
$1)\Rightarrow 2)$ Assume that $N=\oplus _{I}N_{i}$ with
$\left\vert I\right\vert =\aleph ^{\prime }$ and $N_{i}\neq 0$
for every $i\in I$. And let $\left\{ x_{k}\right\} _{K}$ be a
generator set of $N$. For every $k\in K $, there exists a finite
subset $I_{k}\subseteq I$ such that $x_{k}\in \oplus _{i\in
I_{k}}N_{i}$. Thus, $N=\sum_{k\in K}x_{k}R=\sum_{k\in K}\left(
\oplus _{i\in I_{k}}N_{i}\right)=\oplus_{i\in (\cup _{k\in
K}I_{k})}N_{i}$. Therefore, $I=\cup _{k\in K}I_{k}$ and, as each
$I_{k}$ is finite, this means that $\left\vert K\right\vert \geq
\left\vert I\right\vert $.

$2)\Rightarrow 1)$ Assume now that every generator set of $N$ has
cardinality at least $\aleph ^{\prime }$. By Kaplansky's
Theorem \cite[Theorem 26.1]{AF}, $N$ is a direct sum of $\aleph$-generated
submodules, say $N=\oplus _{I}N_{i}$. Let us choose, for every
$i\in I$,
a generator set $%
X_{i}$ of $N_{i}$ of cardinality $\aleph $. Then $\cup _{I}X_{i}$
is a generator set of $N$. Our assumption implies that
$\left\vert \cup _{I}X_{i}\right\vert $ must be at least $\aleph
^{\prime }$. And, as $\aleph <\aleph ^{\prime }$, this means that
$\left\vert I\right\vert \geq \aleph ^{\prime }$
\end{proof}

Recall that a right module $P$ is called hereditary when every
submodule is projective (see e.g. \cite{H}). In particular, every projective module
over a hereditary ring is hereditary.

\begin{corollary}
\label{mingenerators} Let $P$ be a hereditary module.  If
$\Gdim(P)$ is infinite, then it equals the minimum cardinal number
$\aleph $ such that every submodule of $P$ is $\aleph$-generated.
\begin{proof}
As $P$ is hereditary, every submodule of $P$ is projective. Let us
distinguish two possibilities. If $\aleph>\aleph_0$ then the
result is a consequence of the above Lemma.

Assume now that $\aleph=\aleph_0$. Lemma \ref{generators} assures
that $\Gdim(P)\leq \aleph$. Therefore, $\aleph_0\leq \Gdim(P)\leq
\aleph=\aleph_0$.
\end{proof}
\end{corollary}

The following Lemma summarizes some well-known properties of hereditary
modules. We are including a proof for the sake of completeness.

\begin{lemma}\label{hereditary}
Let $P$ be a hereditary right $R$-module. Then
\begin{enumerate}
\item $P$ is non-singular

\item Given an infinite cardinal number $\aleph$, $P$ is
$\aleph$-generated (resp., finitely generated) iff it contains an
essential $\aleph$-generated (resp. finitely generated) submodule.

\item Any submodule of a direct sum of copies of $P$ is a direct
sum of submodules of $P$. In particular, any direct sum of copies
of $P$ is hereditary.
\end{enumerate}
\begin{proof}
\begin{enumerate}
\item Let $x\in P$, $x\neq 0$. Its right annihilator, $r_R(x)$, is
the kernel of the homoeomorphism $f:R\rightarrow xR$ consisting on
right multiplication by $x$. As $xR\subseteq P$, it is projective.
Therefore, $r_R(x)$ is a direct summand of $R_R$ and it cannot be
an essential right ideal. \item Let $N$ be an essential
$\aleph$-generated submodule of $P$. As $P$ is projective, there
exists a splitting epimorphism $\pi:R^{(I)}\rightarrow P$ for some
index set $I$. Let $u:P\rightarrow R^{(I)}$ such that $\pi\circ
u=1_P$. Since $N$ is $\aleph$-generated (resp. finitely
generated), there exists a subset $K\subseteq I$ of cardinality
$\aleph$ (resp. finite) such that $u(N)\subseteq R^{(K)}$. Let
$q:R^{(I)}\rightarrow R^{(K)}$ and $v:R^{(K)}\rightarrow R^{(I)}$
be the canonical projection and injection, respectively. Then
$1_P-(\pi\circ v\circ q\circ u)$ is an endomorphism of $P$ whose
kernel is essential, since it contains $N$. Therefore,
$1_P-(\pi\circ v\circ q\circ u)=0$ (because $P$ is nonsingular)
and this means that $P$ is a direct summand of $R^{(K)}$.

\item This may be proven by following the arguments in \cite[39.7
(2)]{W}. It is interesting, however, to point out that this is
indeed the case even though the definition of hereditary modules
used in \cite{W} is not the one we use in this paper. We are
following the original definition of hereditary modules which goes
back to \cite{H}.
\end{enumerate}
\end{proof}
\end{lemma}
\begin{theorem}\label{dauns-fuchs}
Let $\left\{M_i\right\}_I$ be a family of modules. Then
$\Gdim(\oplus_I M_i)=\sum_I \Gdim(M_i)$.
In particular, for any module
$M$ and any infinite index set $I$,
$\Gdim(M^{(I)})=\max\left\{\Gdim(M),\left\vert I\right\vert\right\}$
\end{theorem}

\begin{proof}
See \cite[Theorem 13]{DF}
\end{proof}

We recall that the cofinality of a cardinal number $\aleph$ is
defined to be the least ordinal number $\alpha$ such that there
exists an injective increasing map $f:\alpha\rightarrow \aleph$
that is cofinal in $\aleph$. I.e., such that for any ordinal
number $\gamma<\aleph$ there exists an ordinal $\beta<\alpha$ with
$f(\beta)\geq\gamma$ (see e.g. \cite[Section 5.4]{R}). The
cofinality of $\aleph$ is always a cardinal number that we will
denote by $\cof(\aleph)$. It is clear that $\cof(\aleph)\leq\aleph$.
A cardinal number $\aleph$ is called regular if
$\cof(\aleph)=\aleph$. Otherwise, $\aleph$ is called singular. An
uncountable cardinal $\aleph$ is said to be (weakly) inaccessible
if it is both regular and limit (i.e., it is not the successor of
any other cardinal). The main result of \cite{DF} shows that
$\Gdim(M)$ is always attained whenever it is not an inaccessible
cardinal.

The following Theorem will be crucial for our first upper bound of
the Goldie dimension of a hereditary module $P$ in terms of its
quasi-injective hull (Theorem \ref{main} below).

\begin{theorem}\label{bound}
Let $P$ be a hereditary module, $M$ a finitely generated submodule
of $P$ and $Q(M)$, the $P$-injective hull of $M$. If $Q(M)$ is
$\aleph$-presented for some infinite cardinal number $\aleph$,
then every submodule of $M$ has a generator set with cardinality
strictly smaller than $\aleph$.

\begin{proof}
Assume on the contrary that there exists a submodule $L$ of $M$
such that every generator set of $L$ has cardinality at least
$\aleph$. We show next that this implies that $L$ contains a
direct sum of $\aleph$ nonzero submodules, say
$\oplus_I L_i$. If $\aleph>\aleph_0$, it is clear by Lemma
\ref{generators}. And, if $\aleph=\aleph_0$, it is a
consequence of the fact that hereditary modules of finite Goldie
dimension are finitely generated by Lemma \ref{hereditary} (2).
Moreover, let us realize that we can assume that each $L_i$ is
finitely generated and therefore, $\oplus_I L_i$ is
$\aleph$-generated. And, adding a complement if necessary (see
\cite[5.21]{AF}), we can also assume that $\oplus_I L_i$ is
essential in $M$.

Let $Q(M)$ be the $P$-injective envelope of $M$. $Q(M)$ is
$\aleph$-presented by assumption. So there exists an epimorphism
$\pi:P^{(A)}\rightarrow Q(M)$ with $\vert A\vert\leq \aleph$ and
$Ker(\pi)$, an $\aleph$-generated module. Using \cite[Theorem
2.2]{GG}, we deduce that there exists a submodule $N$ of $Q$ such
that $\aleph^+$ is attained in $Q/N$. Let $X=\oplus_{j\in J}X_j$
be a direct sum of non-zero modules contained in $Q/N$ with $\vert
J\vert=\aleph^+$. And let $q:Q\rightarrow Q/N$ be the canonical
projection. Then $(q\circ\pi)^{-1}(X)$ is a submodule of $P^{(A)}$
that cannot have a generator set of cardinality at most $\aleph$,
since every generator set of $X$ has clearly cardinality at least
$\aleph^+$. Moreover, $(q\circ\pi)^{-1}(X)$ is projective, since
$P^{(A)}$ is hereditary. Thus, $(q\circ\pi)^{-1}(X)$ is a direct
sum of non-zero modules in cardinality $\aleph^+$, by Lemma
\ref{generators}. Say $(q\circ\pi)^{-1}(X)=\oplus_B Y_b$ with
$\vert B\vert=\aleph^+$.

Let us choose an $\aleph$-generated submodule $N$ of $P^{(A)}$
such that $\pi(N)=\oplus_I L_i$. Then $N+Ker(\pi)$ is essential in
$(q\circ\pi)^{-1}(X)$ because $\oplus_I L_i$ is essential in $M$.
But $N+Ker(\pi)$ is $\aleph$-generated, since so are $N$ and
$Ker(\pi)$. And this means that there exists a subset $B'\subseteq
B$ of cardinality $\aleph<\vert B\vert$ such that
$N+Ker(\pi)\subseteq \oplus_{B'}Y_b$. Let us pick an element $b\in
B\setminus B'$. Then $Y_b\cap\left(N+Ker(\pi)\right)=0$, which is a
contradiction because $N+Ker(\pi)$ is essential in
$(q\circ\pi)^{-1}(X)$.
\end{proof}
\end{theorem}

\begin{corollary}\label{quotients}
Let $\aleph$ be a cardinal number and $P$, an $\aleph$-generated
hereditary module. Let $\aleph'$ be an infinite cardinal number
such that $\cof(\aleph')>\aleph$. If the quasi-injective hull of
$P$ is $\aleph'$-presented, then $\aleph'$ is not attained in
$P/L$ for any submodule $L$ of $P$.

\begin{proof}
Assume on the contrary  that $\aleph'$ is attained in $P/L$ for
some submodule $L$ of $P$. By Lemma \ref{generators} there exists
a submodule of $P/L$ (and thus, a submodule of $P$) satisfying
that every generator set has cardinality at least $\aleph'$.
Reasoning as in the above theorem, we deduce that $\aleph'$ is
also attained in $P$.

Let $\{M_i\}_I$ be the set of finitely generated submodules of
$P$. By Lemma \ref{hereditary}, $P$ is a direct sum of countably
generated submodules of the $M_i$'s, say $P=\oplus_A P_a$, with
$\vert A\vert\leq\aleph<\cof(\aleph')$. Thus, $\aleph'$ is attained
in $P_a$, for some $a\in A$, by Lemma \cite[Lemma 2]{DF}. But
$P_a$ is a submodule of $M_i$ for some $i\in I$. This means that
$\aleph'$ is attained in $M_i$ and thus, $M_i$ contains a
submodule $L$ verifying that any generator set has cardinality at
least $\aleph'$, which is a contradiction with Theorem \ref{bound}.
\end{proof}
\end{corollary}

\begin{corollary}
Let $\aleph$ a cardinal number and $P$ an $\aleph$-generated
hereditary module. Let $\aleph'>\aleph$ be an infinite cardinal
number. If the quasi-injective hull $Q$ of $P$ is
$\aleph'$-presented, then $\Gdim(P/L)\leq \aleph'$ for every
submodule $L$ of $P$.
\begin{proof}
Assume on the contrary that $\Gdim(P/L)>\aleph'$ for some submodule
$L$ of $P$. Then $\Gdim(P)\geq(\aleph')^+$ and thus, $\aleph^+$ is
attained in $P/L$ since it is not an inaccessible cardinal.
Moreover, $\cof((\aleph')^+)=\aleph^+>(\aleph')$. But this
contradicts Corollary \ref{quotients}.
\end{proof}
\end{corollary}

Our next result improves Theorem \ref{bound} when the considered
hereditary module is finitely generated.

\begin{theorem}\label{main}
Let $P$ be a finitely generated hereditary module and $\aleph$, an
infinite cardinal number. If the quasi-injective hull $Q(P)$ of
$P$ is $\aleph$-generated, then $\aleph$ is not attained in $P/L$
for any submodule $L$ of $P$.
\begin{proof}
Assume on the contrary that $\aleph$ is attained in $P/L$ for some
submodule $L$ of $P$. The $\aleph$ is attained in $P$ by the same
reason as in Theorem \ref{bound}. Let $\oplus_I L_i$ be a direct
sum of non-zero submodules of $P$ with $|I|=\aleph$. And let
$\pi:P^{(A)}\rightarrow Q(P)$ be an epimorphism with $|A|\leq
\aleph$. Let $K$ be the kernel of $p$. $K$ is projective and, in
particular, a direct sum of countably generated modules, say
$K=\oplus_J K_j$. Let $\aleph'=\max\{\aleph,\vert J\vert\}$. We
claim that $P$ contains an infinite direct sum of $\aleph'$
non-zero submodules.

By Theorem \ref{dauns-fuchs}, we know that
\[
\Gdim\left(P^{(A)}\right)=\sum_{A} \Gdim(P)=| A|\cdot
\Gdim(P)=\max\left\{\aleph,\Gdim(P)\right\}
\]
and, as we are assuming that $\aleph\leq \Gdim(P)$, we deduce that
$\Gdim(P)=\Gdim\left(P^{(A)}\right)\geq \aleph'$. Let us check that
$\aleph'$ is attained in $P$. If $\aleph'=\aleph$, this is
obvious, since $P$ contains $\oplus_I L_i$ with $\vert
I\vert=\aleph$. So let us assume that $\aleph'=\vert
J\vert>\aleph$. If $\vert J\vert$ is not an inaccessible cardinal,
then $\vert J\vert$ is attained in $P$ by \cite[Theorem 6]{DF}.
Otherwise, $\vert J\vert$ equals its cofinality. And therefore,
$cof\left(\vert J\vert\right)=\vert J\vert>\aleph$ is attained in
$P$, since it is clearly attained in $P^{(A)}$ (see \cite[Lemma
2]{DF}).

On the other hand, $Q(P)$ is $\aleph'$-presented. Therefore,
$\aleph'$ cannot be attained in $P$ by Theorem \ref{bound}. That is a
contradiction which shows that $\aleph$ cannot be attained in
$P/L$.
\end{proof}
\end{theorem}

Next corollary shows that Theorem \ref{main} is particularly
interesting when applied to finitely generated hereditary modules
with countably generated quasi-injective hull.

\begin{corollary}\label{noetherian}
Let $P$ be a finitely generated hereditary $R$-module with
countably generated quasi-injective hull. Then $P$ is a noetherian
module.
\begin{proof}
Let $N$ be any submodule of $P$. The above theorem shows that $N$
must have finite Goldie dimension. Therefore, it is finitely
generated by Lemma \ref{hereditary} (2).
\end{proof}
\end{corollary}

In particular, we get the following corollary for right hereditary
rings that extends results in \cite{GDW,GG1}.

\begin{corollary}
Let $R$ be a right hereditary ring. If the injective envelope of
$R_{R}$ is countably generated, then $R$ is right noetherian.
\end{corollary}

We do not know whether Corollary \ref{noetherian} remains true for
any hereditary right module $P$. However, our next proposition
shows that this is the case when the quasi-injective hull of any
cyclic submodule of $P$ is countably generated.

\begin{proposition}
Let $P$ be a hereditary right $R$-module. If the quasi-injective
hull of any cyclic submodule of $P$ is countably generated, then
$P$ is a countable direct sum of noetherian modules.

\begin{proof}
By Corollary \ref{noetherian}, any cyclic submodule of $P$ is
noetherian. Let $\{N_i\}_{i\in I}$ be the family of all cyclic
submodules of $P$. It is clear that $\oplus_I N_i$ is hereditary
by Lemma \ref{hereditary}(3). Moreover, $P$ is a quotient of
$\oplus_I N_i$ and thus, a direct summand. Finally, $P$ is a
direct sum of submodules of the $M_i$'s, again by Lemma
\ref{hereditary}(3).
\end{proof}
\end{proposition}

\begin{remark}
Let us note that if $P=\oplus_I P_i$ is a direct sum of noetherian
modules, then the Grothendieck category $\sigma[P]$ (see \cite{W}
for the definition) is locally noetherian (i.e., it has a
generator set consisting on noetherian objects). As a module $Q$
in $\sigma[P]$ is injective iff it is a $P$-injective module, we
deduce that the quasi-injective hull of $P$ is
$\Sigma$-quasi-injective. Therefore, the quasi-injective hull of
the hereditary modules considered in the above corollary is
$\Sigma$-quasi-injective.
\end{remark}

We close this Section by proving the following interesting
corollary that gives a positive answer to the question posed by
Dung in \cite{D} and by Dung, G\'omez Pardo and Wisbauer in
\cite{GDW}.

\begin{corollary}\label{GDW}
Let $R$ be a right hereditary ring. If $E(R_{R})$ is finitely
generated, then $R_{R}$ is artinian.
\begin{proof}
We know that $R$ is right noetherian by the above corollary. The
result now follows from \cite[Theorem A]{V}.
\end{proof}
\end{corollary}

\section{Additional results.}

The results given in the above section can only be applied to
hereditary modules. In particular, projective modules over right
hereditary rings. In this section we are going to show that, under
certain additional hypothesis, the given arguments can  be
slightly modified in order to cover other situations. We begin by
showing how to apply them to arbitrary non-singular modules over
right hereditary rings.

\begin{proposition}
Let $R$ be a right hereditary ring, $\aleph$, an infinite cardinal
number and $M$ a non-singular finitely generated right $R$-module.
If the (quasi-)injective hull $Q(M)$ of $M$ is $\aleph$-presented,
then $\aleph$ is not attained in $M$.
\begin{proof}
Let us adapt the arguments given in the proof of Theorem
\ref{bound}. Assume otherwise that $\aleph$ is attained in $M$.
Then $M$ contains a direct sum of non-zero finitely generated
modules, say $\oplus_I M_i$, with $\vert I\vert=\aleph$. We are
assuming that $Q(M)$ is $\aleph$-presented. So there exists an
epimorphism $\pi:R^{(B)}\rightarrow Q(M)$ with $\vert
B\vert=\aleph$ such that $Ker(\pi)$ is $\aleph$-presented. Now all
the arguments given in Theorem \ref{bound} apply to this setup.
\end{proof}
\end{proposition}

The assumption that $M$ is nonsingular in the above proposition is
essential in the proof. The reason is that otherwise we do not
have uniqueness in ($M$-)injective hulls of submodules of $M$ in
$Q(M)$. And therefore, we cannot apply \cite[Theorem 2.2]{GG} in
our arguments (see \cite{O1,O2} for an interesting discussion on
this problem). We are going to show that it is possible to drop
this assumption when we choose a cardinal $\aleph$ with cofinality
$\omega$. Let us note that this situation has a particular
interest since $\aleph_0$ has cofinality $\omega$. First we will
need a technical lemma.

\begin{lemma}\label{ess-cyclic}
Let $p:M\rightarrow N$ be a splitting epimorphism of right
$R$-modules. If $L$ is an essential submodule of $M$, then $p(L)$
is an essential submodule of $N$. In particular, if $M$ contains
an essential cyclic submodule, then so does $N$.
\begin{proof} As $p$ is splitting, there exists a
$u:N\rightarrow M$ such that $p\circ u=1_{M}$. Let $K$ be a
non-zero submodule of $N$. Then $u(K)$ is a nonzero submodule of
$M$. $L\cap u(K)\neq 0$ since $L$ is essential in $M$. Therefore
$0\neq p(L\cap u(K))\subseteq K\cap p(L)$. Thus, $p(L)$ is
essential in $N$.
\end{proof}
\end{lemma}

\begin{theorem}\label{singular}
Let $R$ be a right hereditary ring, $\aleph$, an infinite cardinal
number of cofinality $\omega$ and $M$, a finitely generated right
$R$-module. If the (quasi-)injective hull $Q(M)$ of $M$ is
$\aleph$-presented, then $\Gdim(M)<\aleph$.
\begin{proof}
Assume on the contrary that $\Gdim(M)=\aleph$. As any cardinal
number with cofinality $\omega$ is accessible, $\aleph$ is
attained in $M$. Let $\oplus_A M_{\alpha}\subseteq M$ be a direct
sum of nonzero submodules of $M$ with $\vert A\vert=\aleph'$. And
let $\pi:R^{(I)}\rightarrow Q(M)$ be an epimorphism with $\vert
I\vert\leq \aleph$ and $Ker(\pi)$, $\aleph$-generated.

Let us fix, for any $\alpha\in A$, an (M-)injective hull
$Q_{\alpha}$ of $M_{\alpha}$ within $Q(M)$. As
$\cof(\aleph)=\omega$, Tarski's Lemma \cite[Th\'eor\`eme 7]{T}
assures the existence of a subset $\mathcal{K}\subseteq
\aleph^{\omega}$ with $\vert \mathcal{K}\vert>\aleph$ such that
any $K\in\mathcal{K}$ has cardinality $\aleph_0$ and $K\cap K'$ is
finite if $K\neq K'$. Let $Q_K$ be an ($M$-)injective hull of
$\oplus_K Q_{\alpha}$ in $Q(M)$ for any $K\in \mathcal{K}$. And
let $Z=\sum_{\mathcal{K}}Q_K$.

We claim that $Z$ is not $\aleph$-generated. Otherwise, there
would exist a subset $\mathcal{A}\subseteq {K}$ of cardinality
$\aleph$ such that $Z\subseteq\sum_{K\in \mathcal{A}}Q_{K}$.
Therefore, there would exist an element $K_{0}\in
\mathcal{K}\setminus \mathcal{A}$ such that
$Q_{K_0}\subseteq\sum_{K\in \mathcal{A}}Q_{K}$. As $Q_{K_0}$ is a
direct summand of $Q(M)$, there exists a splitting epimorphism
$q:Q(M)\rightarrow Q_{K_0}$. And $M_{K_0}=q(M)$ is an essential
finitely generated submodule of $Q_{K_0}$ by Lemma
\ref{ess-cyclic}. Therefore, there is a finite subset
$\mathcal{A}'\subseteq \mathcal{A}$ such that
$M_{K_0}\subseteq\sum_{K\in \mathcal{A}'}Q_{K}$.

As $\vert K_0\vert=\aleph_0$ and $K_0\cap K$ is finite for any
$K\neq K_0$, there exist a $k_0\in
K_0\setminus\cup_{\mathcal{A}'}K$. And this means that
$Q_{k_0}\cap \sum_{K\in \mathcal{A}'}Q_{K}\neq 0$, since as
$M_{k_0}$ is essential in $Q_{k_0}$. But this is a contradiction,
because $\oplus_{K\in \mathcal{A}'}\left(\oplus_{k\in
K}Q_k\right)$ is essential in $\sum_{K\in \mathcal{A}'}Q_{K}$ and
$k_0\notin K$ for any $K\in\mathcal{A}'$.

Therefore, we have shown that $Z$ cannot be $\aleph$-generated.
Now, we can use the same arguments as in Theorem \ref{bound} to
get a contradiction.
\end{proof}
\end{theorem}

In particular, we get the following corollary that extends
results in \cite{GDW}.

\begin{corollary}\label{countably}
Let $R$ be a right hereditary ring and $Q$ a countably presented
quasi-injective right $R$-module. Then $Q$ is a countable direct
sum of uniform submodules.
\begin{proof}
Using the arguments of \cite[10.1]{DHSW}, we can write $Q=\oplus_I
Q_{i}$, where $I$ is countable and each $E_i$ is the
quasi-injective hull of a cyclic module. Now, the above theorem
shows that each $Q_i$ has finite Goldie dimension and therefore,
it is a finite direct sum of uniform modules.
\end{proof}
\end{corollary}

\begin{corollary}
Let $R$ be a right hereditary ring and $E$, a finitely presented
injective module. Then every quotient of $E$ is
finite-dimensional.
\begin{proof}
Every quotient of $E$ is finitely generated and injective.
Therefore, the result follows from the above corollary.
\end{proof}
\end{corollary}

\begin{remark}
We recall that a ring $R$ is called right PCI if every cyclic
right module is either free or injective (see \cite{Fa}). It was
proved in \cite{Fa} that any right PCI ring is right hereditary.
Moreover, Damiano \cite{Da} showed that any right PCI ring is
right noetherian (see also \cite{GDW}). The key fact in Damiano's
proof is to show that any finitely presented injective module over
a right hereditary ring has finite Goldie dimension. Therefore,
the above corollary gives an alternative proof of Damiano's
result.
\end{remark}

\section{Final remarks.}
Our arguments show that right hereditary rings with finitely
generated injective hull have finite Goldie dimension and
therefore, are right artinian. Thus solving the question posed in
\cite{Du, GDW}. However we do not know whether they can be extended
to answer the following more general conjecture:

\begin{conjecture}\label{conj}(see e.g. \cite{DHSW}) Let $E$ be a
finitely generated injective module such that every quotient of
$E$ is injective. Then $E$ is a direct sum of uniforms.
\end{conjecture}

Our results show that the answer is "yes" for countably presented
injective modules over right hereditary rings. Furthermore, the
next proposition shows that our arguments can be slightly modified
in order to show that the conjecture is true for the class of
rings having cardinality at most $2^{\aleph_0}$.

\begin{proposition}
Let $R$ be a ring of cardinality at most $2^{\aleph_0}$ and let
$E$ be a countably generated injective module. If every quotient
of $E$ is injective, then $E$ is a direct sum of indecomposable
modules.
\begin{proof}
Assume on the contrary that $E$ is not a direct sum of uniforms.
Then $Gdim(E)\geq \aleph_0$. Let $\oplus_I E_{i}\subset E$ be a
direct sum of nonzero injective submodules with $\vert
I\vert=\aleph_0$. By \cite[Theorem 2.2]{GG} there exists a
submodule $L\subseteq E$ such that $\aleph_1$ is attained in
$E/L$. Let $\oplus_J Q_j$ be a direct sum on non-zero submodules
of $E/L$ with $\vert J\vert=\aleph_1$. As $R$ is right hereditary,
$E/L$ is injective. Therefore, we can choose for any  subset
$X\subseteq J$, an injective envelope $Q_{X}$ of $\oplus_{j\in
X}Q_j$ within $E/L$. Clearly, $Q_X\neq Q_Y$ if $X\neq Y$.

On the other hand $E/L$ is countably generated since it is a
quotient of $E$. And thus, each $E_X$ is also countably generated.
Therefore, $E/L$ is a countably generated module that contains at
least $2^{\aleph_1}$ different countably generated submodules. But
this is a contradiction, since $\vert R\vert\leq 2^{\aleph_0}$.
\end{proof}
\end{proposition}

Let us finish the paper by pointing out that not even for
finitely generated modules over right hereditary rings do we
know the answer to Conjecture
\ref{conj} when $\vert R\vert  >2^{\aleph_0}$. In fact, we do not even
know a counterexample to the
following more general question:

\begin{question}
Let $E$ be a finitely generated module such that any pure quotient
is pure-injective. Is $E$ a direct sum of indecomposable
pure-injective modules?
\end{question}

\noindent {\bf Acknowledgement.} Most of these results were
obtained during the second author's visit to the Center of Ring
Theory and its Applications at Ohio University (CRA). That visit
was supported by the Spanish Ministry of Technology. The author
would like to thank the Ministry for this support and the members
of the CRA for their kind hospitality.


\begin{thebibliography}{99}
\bibitem{AF} F.K. Anderson and K.R. Fuller,
'Rings and Categories of Modules', Springer-Verlag, Berlin (1992).

\bibitem{CR} R.R. Colby and E.A. Rutter, Generalizations
of QF-3 Algebras,
Trans. Amer. Math. Soc. {\bf 153} (1971), 371-386.

\bibitem{Da} R.F. Damiano, A right PCI ring is right
Noetherian, Proc. Amer. Math. Soc. {\bf 77} (1979), 11-14.

\bibitem{D} J. Dauns, 'Modules and
rings', Cambridge University Press, Cambridge (1994)

\bibitem{D1} J. Dauns, Goldie dimension of quotient
modules, J. Australian Math. Soc. {\bf 1} (2001),  no. 1, 11-19.

\bibitem{DF} J. Dauns and L. Fuchs, Infinite Goldie
dimensions, J. Algebra
{\bf 115} (1988), no. 2, 297-302.

\bibitem{DHSW} N.V. Dung, D.V. Huynh, P.F. Smith and R.
Wisbauer, 'Extending modules', Longman, Harlow (1994).

\bibitem{Du} N.V. Dung, A note on hereditary rings or
nonsingular rings with
chain condition, Math. Scand. {\bf 66} (1990), no. 2, 301-306.

\bibitem{ET} P. Erd\"{o}s and A. Tarski, On families of
mutually exclusive
sets, Ann. of Math. {\bf 44} (1943), 315-329.

\bibitem{Fa}C. Faith, When are proper cyclics injectives?,
Pacific J. Math. {\bf 45} (1973), 97-112.

\bibitem{GDW} J.L. G\'{o}mez Pardo, N.V. Dung and R.
Wisbauer, Complete pure-injectivity and endomorphism rings, Proc.
Amer. Math. Soc. {\bf 118} (1993),  no. 4, 1029-1034.

\bibitem{GG1} J.L. G\'{o}mez Pardo and P.A. Guil Asensio,
Endomorphism rings of completely pure-injective modules, Proc.
Amer. Math. Soc. {\bf 124} (1996), 2301--2309.

\bibitem{GG} J.L. G\'{o}mez Pardo and P.A. Guil Asensio,
On the Goldie dimension of injective modules, Proc. Edinburgh
Math. Soc. (2) {\bf 41} (1998), 265-275.

\bibitem{ISS} L. Igl\'{e}sias, C. Santa-Clara, and
F.C. Silva, Unified Min-Max and Interlacing Theorems for Linear
Operators, preprint.

\bibitem{H} D.A. Hill, Endomorphism rings of hereditary
modules, Arch. Math. (Basel) {\bf 28} (1977), no. 1, 45-50.

\bibitem{HJ} K. Hrbaceck and T. Jech, 'Introduction to Set
Theory', Marcel Dekker, New York (1978).

\bibitem{O1}  B.L. Osofsky, A counter-example to a lemma of Skornjakov,
Pacific J. Math. {\bf 15} (1965), 985-987.

\bibitem{O2} B.L. Osofsky, A generalization of
quasi-Frobenius rings, J.
Algebra {\bf 4} (1966), 373-387.

\bibitem{R} Roitman, 'Introduction to
modern set theory', John Wiley \& Sons Inc., New York (1990).

\bibitem{SS} C. Santa-Clara and F.C. Silva, On
infinite Goldie
dimension, J. Algebra {\bf 205} (1998), 617-625.

\bibitem{S} F.L. Sandomierski, Nonsingular rings, Proc.
Amer. Math. Soc. {\bf 19} (1968), 225-230.

\bibitem{V} V. Vinsonhaler, Supplement to the paper
"Orders in QF-3 Rings", J. Algebra {\bf 17} (1971), 149-151.

\bibitem{T} A. Tarski, Sur la decomposition des ensambles
en sous-ensambles
presque disjoints, Fund. Math. {\bf 12} (1928), 188-205.

\bibitem{W} R. Wisbauer, Foundations of Module and Ring
Theory, Gordon and Breach (1991).
\end{thebibliography}
\end{document}